\documentclass[12pt]{amsart}
\usepackage{amssymb,amsmath}

\headsep=1cm \numberwithin{equation}{section}
\newtheorem{theorem}{Theorem}[section]
\newtheorem{lemma}[theorem]{Lemma}

\newtheorem{corollary}[theorem]{Corollary}

\theoremstyle{definition}

\theoremstyle{remark}
\newtheorem{remarks}[theorem]{Remarks}

\newcommand{\Ass}{\operatorname{Ass}}

\newcommand{\grade}{\operatorname{grade}}

\newcommand{\Ext}{\operatorname{Ext}}
\newcommand{\Supp}{\operatorname{Supp}}

\newcommand{\Hom}{\operatorname{Hom}}
\newcommand{\Att}{\operatorname{Att}}

\newcommand{\depth}{\operatorname{depth}}

\newcommand{\fm}{\frak{m}}
\newcommand{\fp}{\frak{p}}
\newcommand{\fq}{\frak{q}}
\newcommand{\fa}{\frak{a}}

\begin{document}
\author[Mafi ]{Amir Mafi  }
\title[Some results on local cohomology modules ]
{Some results on local cohomology modules }

\address{A. Mafi, Arak University, Beheshti St.,
P.O. Box:879, Arak, Iran}
\email{a-mafi@araku.ac.ir}



\begin{abstract} Let $R$ be a commutative Noetherian ring, $\fa$
an ideal of $R$, and let $M$ be a finitely generated $R$-module.
For a non-negative integer $t$, we prove that $H_{\fa}^t(M)$ is
$\fa$-cofinite whenever $H_{\fa}^t(M)$ is Artinian and
$H_{\fa}^i(M)$ is $\fa$-cofinite for all $i<t$. This result, in
particular, characterizes the $\fa$-cofiniteness property of local
cohomology modules of certain regular local rings. Also, we show
that for a local ring $(R,\fm)$, $f-\depth(\fa,M)$ is the least
integer $i$ such that $H_{\fa}^i(M)\ncong H_{\fm}^i(M)$. This
result in conjunction with the first one, yields some interesting
consequences. Finally, we extend the non- vanishing Grothendieck's
Theorem to $\fa$-cofinite modules.
\end{abstract}

\maketitle

\section{Introduction}

Throughout this paper, we assume that $R$ is a commutative
Noetherian ring, $\fa$ an ideal of $R$, and that  $M$ is an
$R$-module. Let $t$ be a non-negative integer. Grothendieck [{\bf
4}] introduced the local cohomology modules $H_{\fa}^t(M)$ of $M$
with respect to $\fa$. He proved their basic properties. For
example, for a finitely generated module $M$, he proved that
$H_{\fm}^{t}(M)$ is Artinian for all $t$, whenever $R$ is local
with maximal ideal $\fm$. In particular, it is shown that
$\Hom_{R}(R/{\fm},H_{\fm}^{t}(M))$ is finitely generated. Later
Grothendieck asked in [{\bf 5}] whether a similar statement is
valid if $\fm$ is replaced by an arbitrary ideal. Hartshorne gave
a counterexample in [{\bf 6}], where he also defined that an
$R$-module $M$ (not necessarily finitely generated) is
$\fa$-cofinite, if $\Supp_{R}(M)\subseteq V(\fa)$ and
$\Ext_R^t(R/\fa,M)$ is a finitely generated $R$-module for all
$t$. He also asked when the local cohomology modules are
$\fa$-cofinite. In this regard, the best known result is that when
either $\fa$ is principal or $R$ is local and $\dim R/{\fa}=1$,
then the modules $H_{\fa}^t(M)$ are $\fa$-cofinite. These results
are proved in [{\bf 8}] and [{\bf 3}], respectively. Melkersson
[{\bf 15}] characterized those Artinian modules which are
$\fa$-cofinite. For a survey of recent developments on
cofiniteness properties of local cohomology, see Melkersson's
interesting article [{\bf 16}]. One of the aim of this note is to
show that, for a finitely generated module $M$, the module
$H_{\fa}^t(M)$ is $\fa$-cofinite whenever the modules
$H_{\fa}^i(M)$ are $\fa$-cofinite for all $i<t$ and $H_{\fa}^t(M)$
is Artinian. This result, in particular, characterizes the
$\fa$-cofiniteness property of local cohomology modules of certain
regular local rings ( see Remark 2.3(ii)). Next, we assume that
$R$ is local  with maximal ideal $\fm$. We prove that $f-\depth
(\fa,M)$, which was introduced in [{\bf 14}], is the least integer
$i$ such that $H_{\fa}^i(M)\ncong H_{\fm}^i(M)$. This result
together with our first mentioned result, in turn yields some
interesting consequences.  Finally, we extend the non-vanishing
Grothendieck's Theorem for $\fa$-cofinite $R$-modules.

\section{The results}

The following theorem describes the behaviour of the cofiniteness
and Artinian property on local cohomology modules.

\begin{theorem} Let $M$ be finitely generated such that
$H_{\fa}^{t}(M)$ is Artinian and that $H_{\fa}^{i}(M)$ is
$\fa$-cofinite for all $i<t$. Then $H_{\fa}^{t}(M)$ is
$\fa$-cofinite.
\end{theorem}

{\bf Proof.} In view of [{\bf 16}, Proposition 4.1], it is enough
to prove that $\Hom_{R}(R/{\fa},H_{\fa}^{t}(M))$ is of finite
length. To prove this, by [{\bf 18}, Theorem 11.38], we consider
the Grothendieck spectral sequence $$
E_{2}^{i,j}=\Ext_{R}^{i}(R/{\fa},H_{\fa}^{j}(M))\underset{i}\Longrightarrow
\Ext_{R}^{i+j}(R/{\fa},M).$$ Since $E_{r}^{0,t}\cong
E_{\infty}^{0,t}$ for $r$ sufficiently large, $E_{\infty}^{0,t}$
is isomorphic to a subquotient of $\Ext_{R}^{t}(R/{\fa},M)$ and,
furthermore, $\ker d_{r-1}^{0,t}\cong E_{\infty}^{0,t}$ for all
$r\geq 3$, where $\ker
d_{r-1}^{0,t}=\ker(E_{r-1}^{0,t}\longrightarrow
E_{r-1}^{r-1,t-r+2})$, we can deduce that $\ker d_{r-1}^{0,t}$ is
finitely generated for $r$ sufficiently large. Next, for all
$r\geq 3$, we have the exact sequence
$$0\longrightarrow \ker d_{r-1}^{0,t}\longrightarrow
E_{r-1}^{0,t}\longrightarrow E_{r-1}^{r-1,t-r+2}.$$ Therefore,
since $E_{r-1}^{r-1,t-r+2}$ is a subquotient of
$E_{2}^{r-1,t-r+2}$, our hypothesis give us that $E_{r-1}^{0,t}$
is finitely generated for $r$ sufficiently large. continuing in
this fashion, we see that $E_{2}^{0,t}$ is finitely
generated; and hence it is of finite length. $\Box$\\

The following corollary is immediate.

\begin{corollary} Let $M$ be finitely generated. Suppose
that the local cohomology module $H_{\fa}^{i}(M)$ is
$\fa$-cofinite for all $i<t$ and that it is Artinian for all
$i\geq t$. Then $H_{\fa}^{i}(M)$ is
$\fa$-cofinite for all $i$.\\
\end{corollary}

\smallskip

\begin{remarks}$(i)$ There is an example in [{\bf 7}, Example 3.4] which shows
that $H_{\fa}^{t}(R)$ is not $\fa$-cofinite for $t=\grade({\fa})$.
However, by the above Theorem, $H_{\fa}^{t}(R)$ is $\fa$-cofinite,
whenever it is Artinian.\\
$(ii)$ Let $(R,{\fm})$ be a regular local ring of characteristic
$p(>0)$ and of dimension $n$. Suppose that $R/{\fa}$ is a
generalized Cohen-Macaulay local ring of dimension $d(>0)$. Then,
by [{\bf 20}, Corollary 1.7] and Theorem 2.1, the local cohomology
modules $H_{\fa}^{i}(R)$ are $\fa$-cofinite if and only if
$H_{\fa}^{n-d}(R)$ is $\fa$-cofinite.\\
\end{remarks}

\smallskip

Let $R$ be a local ring with maximal ideal $\fm$ and let $M$ be a
finitely generated. Following [{\bf 9}], a sequence
$x_{1},\ldots,x_{n}$ of elements of $R$ is said to be an
$M$-filter regular sequence if, for all ${\fp}\in
\Supp(M)\backslash \{\fm\} $, the sequence
$x_{1}/{1},\ldots,x_{n}/{1}$ of elements of $R_{\fp}$ is a poor
$M_{\fp}$-regular sequence. For an ideal $\fa$ of $R$, the
$f-\depth$ of $\fa$ on $M$ is defined as the length of any maximal
$M$-filter regular sequence in $\fa$, denoted by
$f-\depth({\fa},M)$. Here, when a maximal $M$-filter regular
sequence in $\fa$ does not exist, we understand that the length is
$\infty$. For some basic applications of these sequences see [{\bf 2}].\\

\begin{lemma} Let $(R,{\fm})$ be a local ring and suppose that $M$
is finitely generated. Then $f-\depth({\fa},M)=\min\{i\in
\mathbb{N}_{0}: \Supp_{R} H_{\fa}^{i}(M)\nsubseteq \{{\fm \}}\}$.
\end{lemma}

{\bf Proof.} Let $x_{1},\ldots,x_{n}$ be a maximal $M$-filter
regular sequence in $\fa$. If there exists
${\fp}\in\Supp_{R}(H_{\fa}^{i}(M))\backslash \{\fm\}$ for some
$0\leq i\leq {n-1}$, then $x_{1}/1,\ldots,x_{n}/1$ is an
$M_{\fp}$-regular sequence contained in ${\fa}R_{\fp}$. Hence
$H_{\fa}^{i}(M)_{\fp}=0$, which is a contradiction. It therefore
follows that $$ f-\depth({\fa},M)\leq \min\{i\in\mathbb{N}_{0}:
\Supp_{R} H_{\fa}^{i}(M)\nsubseteq \{{\fm}\}\}.$$ Next, by
assumption on $x_{1},\ldots,x_{n}$, there exists ${\fp}\in
\Ass_{R}(M/{(x_{1},\ldots,x_{n})M})\backslash \{{\fm}\}$with
${\fa}\subseteq {\fp}$. Now ${\fp}\in
\Ass_{R}(\Hom_{R}(R/{\fa},M/{(x_{1},\ldots,x_{n})M}))$; and hence\\
${\fp}\in \Ass_{R}(\Ext_{R}^{n}(R/{\fa},M))\backslash \{{\fm}\}$.
Therefore, by [{\bf 11}, Proposition 1.1], ${\fp}\in
\Supp(H_{\fa}^{n}(M))\backslash \{{\fm}\}$, and this completes the
proof. $\Box$\\

\begin{theorem}(see [{\bf 9}, Theorem 3.10] and [{\bf 14}, Theorem 3.1])
 Let $(R,{\fm})$ be a
local ring and suppose that $M$ is finitely generated. Then
$f-\depth({\fa},M)=\min\{i\in \mathbb{N}_{0}: H_{\fa}^{i}(M)\ncong
H_{\fm}^{i}(M)\}$.
\end{theorem}

{\bf Proof.} If $\Supp_{R}(M/{\fa}M)\subseteq \{{\fm}\}$, then
$\sqrt{{\fa}+Ann(M)}={\fm}$; and hence $H_{\fa}^{i}(M)\cong
H_{\fm}^{i}(M)$ for all $i\geq 0$. Therefore
$\min\{i\in\mathbb{N}_{0}: H_{\fa}^{i}(M)\ncong
H_{\fm}^{i}(M)\}=\infty=f-\depth({\fa},M)$; and the result
follows. So, we may assume that $\Supp_{R}(M/{\fa}M)\nsubseteq
\{{\fm}\}$. Let $t=f-\depth({\fa},M)$ and let $x_{1},\ldots,x_{t}$
be an $M$-filter regular sequence in ${\fa}$. Then, by [{\bf 19},
Lemma 1.19], $H_{\fa}^{i}(M)\cong
H_{(x_{1},\ldots,x_{t})}^{i}(M)\cong H_{\fm}^{i}(M)$, for all
$i<t$. On the other hand, by Lemma 2.4, the $R$-module
$H_{\fa}^{t}(M)$ is not isomorphic with $H_{\fm}^{t}(M)$. It
therefore follows, by [{\bf 9}, Theorem 3.10]. $\Box$\\

\begin{remarks} Let $M$ be finitely generated. Then\\
$(i)$ in view of Theorem 2.1 and Theorem 2.5, it is clear that if
$(R,{\fm})$ is a local ring, then $H_{\fa}^{i}(M)$ is
$\fa$-cofinite for all $i$ less than $f-\depth({\fa},M)$;\\
$(ii)$ it follows immediately from [{\bf 9}, Theorem 3.10] and
Theorem 2.5 that if $(R,{\fm})$ is local and $H_{\fa}^{i}(M)$ is
Artinian for all $i<t$, then $H_{\fa}^{i}(M)\cong H_{\fm}^{i}(M)$
for all $i<t$.\\
\end{remarks}

\smallskip

The following lemma is needed in the proof of the next theorem.
Note that if we replace $\fa$ by the zero ideal in the lemma, then
the
Grothendieck's Theorem [{\bf 4}, p.88] immediately follows.\\

\begin{lemma} Let $M$ be $\fa$-cofinite. Then for every
maximal ideal $\fm$ of $R$ and for all $t$, $H_{\fm}^{t}(M)$ is
Artinian.
\end{lemma}

{\bf Proof.} Since $H_{\fm}^{t}(M)$ is an $\fa$-torsion module, by
[{\bf 13}, Theorem 1.3], it is enough to prove
$0:_{H_{\fm}^{t}(M)}\fa$ is Artinian. Let $\Phi(-)$ denote the
composite functor $\Hom_{R}(R/{\fa},H_{\fm}^0(-))$. We get a
spectral sequence arising from the composite functor as:

$$E_{2}^{i,j}=\Ext_{R}^{i}(R/{\fa},H_{\fm}^{j}(M))\Longrightarrow
(R^{i+j}\Phi)(M).$$

Now, we use induction on $j$ (with $0\leq j\leq t$) to show that
$E_{2}^{0,t}$ is Artinian. Let $0\leq j<t$ and suppose that the
result has been proved for smaller values of $j$.( Note that the
case $j=0$ was proved in [{\bf 15}, Corollary 1.8].) We can apply
[{\bf 15}, Theorem 1.9] and use a similar argument as in the proof
of Theorem 2.1, to see that $\ker d_{r-1}^{0,j+1}$ is Artinian for
$r$ sufficiently large. On the other hand, by induction,
$E_{r-1}^{r-1,j-r+3}$ is Artinian. It now follows that
$E_{2}^{0,j+1}$ is Artinian. This complete the inductive step. In
particular $E_{2}^{0,t}$ is Artinian. $\Box$

\smallskip

In the next result, we will use the concept of attached prime
ideals. For more details in this subject the reader is referred to
[{\bf 10}] or the appendix
to $\S ${\bf 6} in [{\bf 12}].\\

\begin{theorem}  Let $(R,{\fm})$ be a local ring and let $M$ be a
module of dimension $d$. If $H_{\fm}^{d}(M)$ is an Artinian
module, then if $\fp$ is any of its attached prime ideals, one has
$\dim R/{\fp}\geq d$.
\end{theorem}

{\bf Proof.} From the right exactness of $H_{\fm}^{d}(-)$ on
modules of dimension $\leq d$, we get $H_{\fm}^{d}(M/{\fp}M)\cong
H_{\fm}^{d}(M)/{{\fp}H_{\fm}^{d}(M)}$, which is $\neq 0$, since
${\fp}$ is an attached prime ideal of $H_{\fm}^{d}(M)$. But
$M/{\fp}M$ is a module over $R/{\fp}$. Therefore $\dim R/{\fp}\geq
d$. $\Box$

\smallskip

In the following theorem, which establishes the non-vanishing
Grothendeick Theorem for $\fa$-cofinite modules.

\begin{theorem} Let $(R,{\fm})$ be a local ring and let $M$ be a
non-zero $\fa$-cofinite $R$-module of dimension $n$. Then
$H_{\fm}^{n}(M)\neq 0$.
\end{theorem}

{\bf Proof.} Firstly note that, in view of the hypotheses,
$0:_{M}\fa$ is a finitely generated $R$-module of dimension $n$.
Now, we prove the theorem by induction on $n(\geq 0)$. If $n=0$,
then $0:_{M}\fa$ is Artinian; and hence, by [{\bf 13}, Theorem
1.3], $M$ is Artinian. Therefore $H_{\fm}^{0}(M)=M\neq 0$.\\
Suppose, inductively, that $n\geq 1$ and the result has been
proved for $n-1$. We may assume that $M$ is $\fm$-torsion free.
Also, by [{\bf 15}, Corollary 1.4], we may assume that $\Ass(M)$
is a finite set. Then, there exists a non-zero divisor $x\in
{\fm}$ on $M$. Suppose the contrary that $H_{\fm}^{n}(M)=0$. Then,
for any such $x$, we can consider the exact sequence
$0\longrightarrow M\overset{x}\longrightarrow M\longrightarrow
M/{{x}M}\longrightarrow 0$ to see that
$H_{\fm}^{n-1}(M)/{{x}H_{\fm}^{n-1}(M)}\cong
H_{\fm}^{n-1}(M/{{x}M})$,\\ $n-1=\dim
(0:_{M}\fa)/{{x}(0:_{M}\fa)}\leq \dim (0:_{M/{{x}M}}\fa)=\dim
M/{{x}M}\leq n-1$, and that, by [{\bf 15}, Remark(a)], $M/{{x}M}$
is $\fa$-cofinite. Therefore, by induction hypothesis,
$H_{\fm}^{n-1}(M)/{{x}H_{\fm}^{n-1}(M)}\neq 0$. Note that, by
Lemma 2.7, $H_{\fm}^{n-1}(M)$ is Artinian. If ${\fm}\notin \Att
H_{\fm}^{n-1}(M)$, then, for any $$y\in{\fm}\setminus
\bigcup_{{\fp}\in \Att H_{\fm}^{n-1}(M)}{\fp}\bigcup
\bigcup_{{\fq}\in\Ass(M)}{\fq},$$ we have
$H_{\fm}^{n-1}(M)={y}H_{\fm}^{n-1}(M)$, which is a contradiction.
Thus $\fm \in \Att H_{\fm}^{n-1}(M)$. Let $\Att
H_{\fm}^{n-1}(M)=\{{\fp}_{1},\ldots,{\fp}_{t},{\fm}\}$ and let
$z\in
{\fm}\setminus\bigcup_{i=1}^{t}{\fp}_{i}\bigcup\bigcup_{{\fq}\in\Ass(M)}{\fq}$.
Then, by the above argument, we have
$H_{\fm}^{n-1}(M)/{{z}H_{\fm}^{n-1}(M)}\cong
H_{\fm}^{n-1}(M/{{z}M})$. Hence, by [{\bf 17}, Proposition 5.2],
$\Att H_{\fm}^{n-1}(M/{{z}M})=\Supp(R/{(zR)})\cap\Att
H_{\fm}^{n-1}(M)=\{{\fm}\}$. Therefore, by [{\bf 1}, Corollary
7.2.12], $H_{\fm}^{n-1}(M/{{z}M})$ has finite length. If we show
that $H_{\fm}^{n-1}(M/{{z}M})=0$, then we achieved at the required
contradiction. To this end, first let $n=1$. Then we have the
exact sequence $$0\rightarrow H_{\fm}^{0}(M)\overset{z}\rightarrow
H_{\fm}^{0}(M)\rightarrow H_{\fm}^{0}(M/zM)\rightarrow
H_{\fm}^{1}(M).$$ By our hypothesis
$H_{\fm}^{0}(M)=0=H_{\fm}^{1}(M)$; and so $H_{\fm}^{0}(M/zM)=0$.
Now, we assume that $n>1$. Then, Theorem 2.8 implies that attached
prime ideals of $ H_{\fm}^{n-1}(M/zM)$ is empty; and so
$H_{\fm}^{n-1}(M/zM)=0$. $\Box$\\

{\bf Acknowledgment}: The author is deeply grateful to the referee
for his or her careful reading of the manuscript. The author also
thank the referee for proposing Theorem 2.8.


\end{document}